\input amstex
\voffset=-0.6in
\document
\documentstyle{amsppt}
\magnification1200
\TagsOnRight
\endsubjclass
\topmatter

\title{The concentration index of subharmonic functions of infinite order}
\endtitle
\author{
 Markiyan Hirnyk\footnote"*"{Hirnyk=Girnyk} }
\endauthor

\abstract{The purpose of this paper is to introduce into
consideration an analogue of the concentration index in the class of
subharmonic functions of infinite order. The one in the case of
finite order is used in the
interpolation theory.}
\endabstract
\rightheadtext{Concentration index of subharmonic functions of
infinite
 order}
\endtopmatter

\redefine\delta{\varkappa}
We  use the standard notation of the potential theory and the value
 distribution theory \cite{1, 2}, nevertheless we recall some of them. We denote by
 $\mu_u $ the Riesz measure of a subharmonic function $u$. We put
 $C(z,t)=\{w:|w-z| \le t\}$,\, $n(z,t)=\mu_u(C(z,t))$,\,$\,n(r)=n(0,r)$, and
 $B(r,u)$ the maximum of $u $ on the disk $C(0,r)$.
Without loss of generality we may assume $u(0)=0$ and $n(1)=0$. The set of
 all the subsets of $[1,\infty)$, having finite logarithmic measure, is
 denoted by $FLM$:if $S\in FLM,$ then $ \int_1^{\infty}\chi_S(t)\,d\log (t)< \infty$, where $\chi_S $ is the characteristic function of $S$ . We denote by $M$ positive constants.

  The concentration index of an entire function of finite order was
  introduced into consideration implicitly by Levin \cite{3} and explicitly
 by Krasichkov \cite{4}, who studied its properties. The specific
case of zero order was considered in \cite{5,6}.

 We define the concentration index $I(z,u)$ of a subharmonic function of
 infinite order by the formula
$$
I(z,u)=-\int\limits_0^{|z|/ \log^{\delta}n(|z|)} n(z,t)/t \,dt,
$$
where a real number $\delta>0$.

We prove
\proclaim{Theorem} Let $u$ be a subharmonic function of infinite order, $r=|z|
$. Then
$$
u(z)=I(z,u)+\exp\left(o(N(r)))+O(B(r,u)\right), z\to \infty, r\notin S \in FLM.   \tag1
$$
and
$$
I(z,u)=o\left(\exp (o(N(r)))\right),z \to \infty, z \notin E,                      \tag2
$$
which is such that for every $r\notin S \in FLM$ there exists an at most countable
set of disks $C(z_j,r_j)$, having the following properties:
$$
\cup_jC(z_j,r_j)\supset E\cap \{w:r<|w|<r+\Delta\},                    \tag3
$$
and
$$
\sum_{|z_j|\in [r,r+\Delta]} r_j=o(\Delta),\,r\to \infty,\, r\notin S\in FLM, \tag4
$$
where $\Delta=r/\log^{\delta}n(r)$.
\endproclaim

{\it Proof} . We start with the construction of a subharmonic function $v$ such that
the Riesz measure $\mu_v=\mu_u$ and the growth of the function $v$ is minimal
in some sence.

Let real numbers $\delta $ and $\eta$ satisfy the inequalities $ 0<\delta<
\eta$. Following \cite{7}, we put
$$
v(z)=\int_{\Bbb C}\log |E(z/ \xi,[\log^{1+\eta}n(|\xi|)])|\, d\mu_u(\xi),
\tag4
$$
where $E(z,p)$ is the Weierstrass primary factor of genus $p$. The integral in
the right-hand side of (4) converges uniformly on every compact subset of $ {\Bbb C}$.
This known statement will be proved below too.

We represent $v(z)$ as the sum
$$
v(z)=v_1(z)+v_2(z)+v_3(z)+v_4(z)+v_5(z),  \tag5
$$
where $(R=r+\Delta)$
$$
v_1(z)=\int_{C(z,\Delta)}\log|1-z/\xi|\,d\mu_u(\xi),
$$
$$
v_2(z)=\int_{C(z,\Delta)} \Re \sum \limits _{j=1}^{j=[\log^{1+\eta}n(r)]} j^{-1}(z/\xi)^j \,d\mu_u(\xi),
$$
$$
 v_3(z)=\int_{C(0,r)\setminus C(z,\Delta)} \log|E(z/\xi,[\log^{1+\eta}n(|\xi|)])\,d\mu_u(\xi),
$$
$$
v_4(z)=\int_{C(0,R)\setminus (C(0,r) \setminus C(z,\Delta))} \log|E(z/\xi,[\log^{1+\eta}n(|\xi|)])\,d\mu_u(\xi),
$$
$$
v_5(z)=\int_{\Bbb C \setminus C(0,R)}  \log|E(z/\xi,[\log^{1+\eta}n(|\xi|)])\,d\mu_u(\xi).
$$

We first prove two estimates we will need later on. Applying the
Borel-Nevanlin\-na Theorem \cite{2,\,p.120} with $u(r)= \log \log
n(\exp(r)),\,\, \varphi (u)= \exp(-\delta u + \log M)$ we obtain
$$
n\left(r\left(1+\frac{M}{\log ^{\delta}n(r)}\right)\right)\leq n(r)^e, \,\, r\notin S \in FLM.  \tag6
$$

By Lemma 1.1 \cite{2,p.433} with $\varepsilon=1,\varphi (t)=N(\exp(t))$, we have
$$
n(r)\leq N(r)^2,\,\, r\notin S\in FLM. \tag7
$$

We now turn to the estimation $v_1(z)$. We denote by $\nu(z,t)$ the
measure $$\mu_u(C(0,t) \cap C(z,t)). $$ Providing $u(z)>-\infty$, we
have
$$\split
v_1(z)=\int\limits_0^{\Delta}\log t\,dn(z,t)+\int\limits_0^R\log \frac{1}{t}
\,d\nu(z,t)=-I(z,t)+\log\Delta\, n(z,t)\\+\int\limits_0^R \frac{\nu(z,t)}{t}
\,dt+\log\frac{1}{R}\,\nu(z,R).\endsplit\tag8
$$
Above we wrote $v_1(z)$ as the sum of the Stieltjes integrals and integrated
by parts. Applying (6), (7), we obtain
$$
\split
|\log \Delta n(z,\Delta)|\leq n(R)(\log r +\delta \log \log n(r))=\\=
O(n(r)^{e+1}\log r)=O(N(r)^{2e+3}),\,r\to \infty,\, r\notin S\in FLM.\endsplit
\tag9
$$
Likewise,
$$
\split
\left|\log \frac{1}{R}\nu(z,R)\right|\leq n(R)\log R =O(n(r)^e\log r)=\\=
O(N(r)^{2e+1}),\,r\to \infty,\, r\notin S\in FLM.\endsplit
\tag10
$$
Next,
$$
\split
\int\limits_0^R\frac{\nu(z,t)}{t}\,dt\leq n(R)\log R \\=O(N(r)^{2e+1}),\,r\to \infty,\, r\notin S\in FLM.\endsplit
\tag11
$$
Combining (8)-(11), we obtain
$$
v_1(z)=I(z,u)+O(N(r)^{2e+3}),\,r\to \infty,\, r\notin S\in FLM.
\tag12
$$

We will need the elementary inequality
$$
\sum\limits_{j=1}^{j=p}j^{-1}|w|^j\leq a^p(2+\log p),\tag13$$ which
holds under the assumptions $|w|<a$ and $a>1$.

Applying (13) to the estimation $v_2$, we have
$$\split
|v_2(z)|\leq 2\left(\frac{R}{r-\Delta}\right)^{log^{1+\eta}n(R)}\log( \log n(R))\,n(R)\leq \\
\leq 2\left(1+3 \log^{-\delta}n(r)^{\log ^{1+\eta}n(R)}\right)=O(n(r)^{e+1}\exp(O(1)\log^
{1+\eta-\delta}n(r)))=\\=O(\exp(o((N(r))),\, r\to \infty,\,r\notin S\in FLM.\endsplit
\tag14
$$
We now take up the consideration of $v_3$. In view of the inequality
$$\multline
|\log|1-z/\xi||\leq \max\left(\left|\log \frac \Delta
r\right|,\log(1+r)\right)\leq 2(\log \log n(r)+\log r),\endmultline
\tag15
$$
which holds on the set $C(0,r)\setminus (C(z,\Delta)\cup C(0,1))$, and (13),
we obtain
$$
\multline
|v_3(z)|\leq\\ \leq n(r)(2\log \log n(r))+2\log r+(2\log \log n(r)+2)\left(\frac {r}{r_0}
\right)^{(\log n(r_0))^{1+\eta}},\endmultline
$$
where
$$
\left(\frac{r}{r_0}\right)^{(\log n(r_0)^{1+\eta})}=\max_{1\leq |\xi|\leq r}\left(\frac {r}{|\xi|}\right)
^{  \log ^{1+\eta}n(|\xi|)}), \,1\leq r_0\leq r,
$$
and $r_0$ is the greatest such number. It exists, because the function $n(r)$
is upper semicontinuous on $[1,r]$, as a nondecreasing and right-continuous
function. We easily see $r_0 \to \infty$ as $r\to \infty$. Taking into account
the inequality
$$
N(r)\geq \int\limits_{r_0}^r \frac{n(t)}{t}\,dt \geq n(r_0)\log \frac {r}{r_0}
$$
and (7), we have
$$\multline
|v_3(z)|\leq 4n(r)(\log \log n(r)+\log r)\exp\left(\frac{\log^{1+\eta}n(r_0)}{n(r_0)}
N(r)\right)=\exp (o(N(r))),\\r\to \infty, r\notin S\in FLM.\endmultline \tag16
$$

 The next term $v_4(z)$ is estimated somewhat in another way. If $\xi \in
C(0,R)\setminus(C(0,r)\cup C(z,\Delta))$, then
$$
\left|\log\left|1-\frac {z}{\xi}\right|\right| \leq \left|\log \frac {\Delta}{R}\right|\leq \log \log n(r).\tag17
$$
From (6), (13), and (17) we conclude (compare with(14))
$$\multline
|v_4(z)|\leq n(R)\left(\log \log n(r)+\left(\frac R r\right)\right)^{\log^{1+\eta}n(R)}(\log \log n(R)
+2))\leq \\ \leq n(r)^e \left(\log \log n(r)+\left(1+\frac{1}{\log ^{\delta}n(r)}\right)^
{(e\log n(r))^{1+\eta}}\right)(\log \log n(r)+3))=\\=O(\exp(o(N(r)))),\\
r\to \infty, r\notin S\in FLM.\endmultline\tag18
$$

Finally, we estimate $v_5(z)$. Applying the inequality \cite{3, p.21}
$$
|\log|E(w,p)||\leq |w|^{p+1},\,\text {when}\,\, |w|\leq\frac{p}{p+1},
$$
we obtain
$$\multline
|v_5(z)|\leq \int\limits_{\Bbb C \setminus
C(0,R)}\left(\frac{r}{|\xi|}\right)^
{\log^{1+\eta}n(|\xi|)+1}\,d\mu_u(\xi)\leq \\ \leq \int\limits_{\Bbb
C \setminus C(0,R)}\left(\frac{r}{R}\right)^{\log^{1+\eta}
n(|\xi|)+1}\,d\mu_u(\xi)\leq\int\limits_{\Bbb C \setminus
C(0,R)}\left(\frac{r}{R}\right)^ {\log^{\delta}n(r)\log
^{1+\eta-\delta}n(|\xi|)}\,d\mu_u(\xi) \\\leq \int\limits_{\Bbb C
\setminus C(0,R)}
2^{-\log^{1+\eta-\delta}n(|\xi|)}\,d\mu_u(\xi)=\int\limits_{R}^{\infty}
2^{-(\log n(t))^{1+\eta-\delta}}\,dn(t)=O(1), \,r\to\infty.
\endmultline\tag19
$$

Combining (12),(14),(16),(18), and (19), we have
$$
v(z)=I(z,u)+\exp(o(N(r))),\,r\to \infty,r\notin S\in FLM,\tag20
$$
i.e. the modulus of the difference $v(z)-I(z,u)$ is bounded by a
nondecreasing function $V$, which is such that $V(r)=\exp
(o(N(r))),\,r\to\infty,\, r\notin S\in FLM$ and $N(r)^{2e}=o(V(r)),
\,r\to \infty$ (We can increase $V(r)$ in need).

The next step consists in the proof of claims (1)-(4). We will use a method by
Hayman \cite{5}. A point $z$ is said to be $(\beta, s)$-light with respect to
a measure $\mu$ if for every $t\in (0,s)$ the inequality $n(z,t)< \beta t$ holds.
We denote by $LP(\beta, s, \mu)$ the set of such points. We put $s(z)=s(|z|)=r/\log^{\delta}n(r)$, $\mu=\mu_u$.
  We choose $\beta (z)=\beta (|z|)$ in such a way that
  $$
  \beta (z)s(z)=o(V(r)),\,r\to \infty,\,r\notin S \in FLM , \tag21
  $$
  $$
  6N(r)^{2e}=o(\beta(r)s(r))  ,\,r\to \infty,\,r\notin S \in FLM , \tag22
  $$

  For instance, we can put $\beta(r)s(r)=(V(r)N(r)^{2e})^{1/2}.$
 If a point
 $z\in LP(\beta, s, \mu)$,\newline then, applying (21), we have
$$\multline
|I(z,\mu)|=\int\limits_0^{\Delta}\frac{n(z,t)}{t}\,dt\leq\int\limits_0^{\Delta}
\beta \,dt=\\=\beta(z)s(z)=o(\exp(o(N(r)))),
r\to \infty,\,r\notin S\in FLM.\endmultline
$$

If a point $z$ is heavy (i. e. $z\in HP(\beta, s, \mu)= \Bbb C \setminus LP(\beta, s, \mu)$)
, then there exists a real number $r_z\in (0, s)$, such that $n(z,r_z)\geq \beta r_z$.
We obtain a cover $\{C(z,r_z)\}$ of the set $HP(\beta, s, \mu)$. Applying the
Besicovitch-Landkof Theorem \cite{9, p.246}, we can choose an at most countable
subcover $\{C(z_j, r_j)\}$ of multiplicity less than or equal to 6.

We note that if $t\in[r,R]$, then
$$
\frac {r}{\log^{\delta}n(R)}\leq s(t)\leq \frac {R}{\log^{\delta}n(r)},
$$
and thus $s(t)\sim \Delta, \,r\to \infty,\,r\notin S\in FLM.$

Because of this we have
$$
\multline
\sum\limits_{|z_j|\in [r,R]}n(z_j,r_j)\leq 6n(R+\Delta)=6n\left(r+\frac{2r}{\log^{\delta}n(r)}\right)\\
\leq 6n(r)^e \leq 6N(r)^{2e},\,r\notin S\in FLM.
\endmultline \tag23
$$

On the other hand,
$$
\sum\limits_{|z_j|\in [r,R]}n(z_j,r_j)\geq\sum\limits_{|z_j|\in [r,R]}
\beta(|z_j|)r_j \geq \beta(r)\sum\limits_{|z_j|\in [r,R]}r_j.\tag24
$$
Comparing (23), (24) and also using (22), we obtain
$$
\multline\sum\limits_{|z_j|\in [r,R]}r_j\leq
6N(r)^{2e}\beta(r)^{-1}= o(\Delta),\,r\to \infty,\,r\notin S\in
FLM.\endmultline
$$
To complete the proof of the theorem, we show
$$
|u(z)-v(z)|\leq M\left(B(r,u)+\exp(o(N(r)))\right)\,r\to \infty,\,r\notin S\in FLM.
\tag25
$$

As Goldberg proved in \cite{7}(He considered only the case of entire functions,
 but his result and  naturally changed proof are true for subharmonic functions too.),
$$
B(r,v)\leq\exp(o(N(r))),\,r\to \infty,\,r\notin S\in FLM.
$$
Combining this and Theorem 4.4\cite{10}, we obtain
$$\multline
u(z)-v(z)\leq M\,T(r,u-v)\leq M\,(T(r,u)+T(r,-v))=\\=M(T(r,u)+T(r,v)\leq
M(T(r,u)+\exp(o(N(r)))),\,r\to \infty,\,r\notin S\in FLM.\endmultline
$$
Above we used the First Main Theorem of the value distribution
theory  . It should be noted we  have no exceptional set of disks,
because the function $u-v$ is harmonic. We can apply the same
arguments to $v-u$ too, thus we prove (25).

I am grateful to Professor M. Zabolotskii for the statement of the prob\-lem
and the par\-ti\-cipants of Lviv se\-minar on comp\-lex ana\-lysis for useful discussions.

My special thanks to the referee for a careful review of the paper.

\centerline {References}

\item{1.}  W.K.~Hayman and P.B.~Kennedy, Subharmonic Functions. Vol.1.
Academic Press. London-New York-San Francisco, 1976.
\item{2.}A.A.~Goldberg and I.V. Ostrovskii, Value Distribution of
Meromorphic Functions.  Nauka, Mos\-cow, 1970 (Russian).
\item{3.}B.Ya.Levin, Distribution of Zeros of Entire Functions, GITTL, Mosc
ow, 1956 (Russian).
\item{4.}I.F. Krasichkov, Lower estimates for entire functions of finite order.
Sib. Math. Zh., v. 6 (1965) , N 4, 840-861 (Russian).
\item{5.} A.A.Goldberg and N.V. Zabolotskij, The concentration index of a
subharmonic function of zero order, Mat. Zametki, v. 34 (1983), N 2, 227-236 (Russian).
\item{6.} N.V. Zabolotskij and S.Yu. Favorov, Asymptotic formulae for
subharmonic in ${\Bbb R}^m$ functions of zero order, Teor. Funfkts., Funkts.
Anal. Ikh Prilozh., v.47  (1987), 125-128 (Russian).
\item{7.}A.A. Goldberg, The representation of a meromorphic function as a
quotient of entire functions, Izvestiya Vuzov, (1972) N 10, 13-17 (Russian).
\item{8.}W.K. Hayman, Questions of regularity connected with the Phragm\'en-
Lindel\"of principle, J. Math. Pures et Appl., v. 32 (1956), N 2, 115-126.
\item{9.} N.S. Landkof, Foundations of Modern Potential Theory, Nauka, Moskow
, 1966 (Russian).
\item {10.} V. Ya. Eiderman, Estimates of potentials and $\delta$-subharmonic
 functions outside exceptional sets, Izvestiya RAN , v. 61 (1997), N 6, 181-218
 (Russian).
\end